\documentclass{article}

\usepackage{arxiv}
\usepackage{amsmath, amsthm, amssymb, epsfig, sectsty, graphicx, float, url}

\usepackage{algorithm, algorithmicx, algpseudocode}

\usepackage{float}
 \usepackage{subcaption}
 \usepackage{caption}

\usepackage[utf8]{inputenc} % allow utf-8 input
\usepackage[T1]{fontenc}    % use 8-bit T1 fonts
\usepackage{url}            % simple URL typesetting
\usepackage{booktabs}       % professional-quality tables
\usepackage{amsfonts}       % blackboard math symbols
\usepackage{nicefrac}       % compact symbols for 1/2, etc.
\usepackage{microtype}      % microtypography
\usepackage{lipsum}
\numberwithin{equation}{section}

\usepackage{amsmath, amsthm, amssymb, epsfig, sectsty, graphicx, float, url}
\usepackage{algorithm, algorithmicx, algpseudocode}
\usepackage{float}
 \usepackage{subcaption}
 \usepackage{caption}
 
\usepackage[colorinlistoftodos,bordercolor=orange,backgroundcolor=orange!20,line
color=orange,textsize=scriptsize]{todonotes}
\catcode`\@=11
\catcode`\@=12

\numberwithin{equation}{section}

\usepackage{bm}

\theoremstyle{plain}
\newtheorem{theorem}{Theorem}[section]
\newtheorem{lemma}[theorem]{Lemma}

\newtheorem{proposition}[theorem]{Proposition}

\theoremstyle{definition}
\newtheorem{definition}[theorem]{Definition}

\newcommand{\R}{\mathbb{R}}

\newcommand{\A}{\mathbb{A}}

\newcommand{\N}{\mathbb{N}}

\newcommand{\X}{\mathbb{S}}
\newcommand{\XX}{\mathcal{X}}
\newcommand{\VV}{\mathcal{V}}
\newcommand{\FF}{\mathcal{F}}

\newcommand{\E}{\mathbb{E}}

\usepackage[utf8]{inputenc} % allow utf-8 input
\usepackage[T1]{fontenc}    % use 8-bit T1 fonts

\usepackage{url}            % simple URL typesetting
\usepackage{booktabs}       % professional-quality tables
\usepackage{amsfonts}       % blackboard math symbols
\usepackage{nicefrac}       % compact symbols for 1/2, etc.
\usepackage{microtype}      % microtypography

\usepackage{hyperref}
\hypersetup{
    colorlinks=true,
    linkcolor=blue,
    citecolor=red}

\title{From Convex Optimization to MDPs: A Review of First-Order, Second-Order and Quasi-Newton Methods for MDPs}

\author{%
   Julien Grand-Cl{\'e}ment \\
ISOM Department, HEC Paris \\
   \texttt{grand-clement@hec.edu}
}

\begin{document}
\maketitle
\vspace{4mm}

\begin{abstract}
In this paper we present a review of the connections between classical algorithms for solving Markov Decision Processes (MDPs) and classical gradient-based algorithms in convex optimization. Some of these connections date as far back as the 1980s, but they have gained momentum in recent years and have lead to faster algorithms for solving MDPs.  In particular,  two of the most popular methods for solving MDPs, Value Iteration and Policy Iteration, can be linked to first-order and second-order methods in convex optimization.  In addition, recent results in quasi-Newton methods lead to novel algorithms for MDPs, such as Anderson acceleration.  By explicitly classifying algorithms for MDPs as first-order, second-order, and quasi-Newton methods, we hope to provide a better understanding of these algorithms, and, further expanding this analogy,  to help to develop novel algorithms for MDPs, based on recent advances in convex optimization.
\end{abstract}

%\tableofcontents


\section{Introduction}

Markov Decision Process (MDP) is a common framework modeling dynamic optimization problems, with applications ranging from reinforcement learning \citep{mnih2015human} to healthcare \citep{Goh,grand2020robust} and wireless sensor networks \citep{alsheikh2015Markov}.  Most of the algorithms for computing an optimal control policy are variants of  two algorithms: \textit{Value Iteration} (VI) and \textit{Policy Iteration} (PI).  Over the last 40 years,  a number of works have highlighted the strong connections between these algorithms and methods from \textit{convex optimization}, even though computing an optimal policy is a non-convex problem.  Most algorithms in convex optimization can naturally be classified as \textit{first-order}, \textit{second-order} and \textit{quasi-Newton} methods, if the iterates rely on gradients and/or Hessian computations \citep{boyd-2004}.   The goal of this  paper is to outline a unifying framework for the classification of algorithms for solving MDPs.  In particular, we present a systematic review of the connections between Value Iteration and first-order methods,  between Policy Iteration and second-order methods, and between variants of Value Iteration and quasi-Newton methods.  We hope that this unifying view can help to develop novel fast algorithms for solving MDPs, by extending the latest advances in convex optimization to the MDP framework. On the other hand, solving MDPs through the lens of convex optimization motivates novel interesting questions and challenges in optimization, as the operators and objective functions do not satisfy classical structural properties (e.g, convexity and/or differentiability).

\paragraph{Outline.} We introduce the MDP framework as well as the classical Value Iteration and Policy Iteration algorithms in Section \ref{sec:MDP}.  We highlight the recent connections between Value Iteration and first-order methods (Gradient Descent) in Section \ref{sec:VI-GD}. The relations between Policy Iteration and second-order methods (Newton's method) are presented in Section \ref{sec:PI-newton}. We review Anderson Value Iteration, a quasi-Newton methods for MDPs, in Section \ref{sec:anderson-mdp}. For the sake of completeness, in Appendix \ref{sec:convex-opt}, we present a detailed review of the results for Gradient Descent (along with acceleration and momentum),  Newton's method and quasi-Newton methods in convex optimization. 

\paragraph{Notations.} In this paper, $n$ and $A$ denote integers in $\N$.  The notation $\Delta(A)$ refers to the simplex of size $A$.  We write $[n]$ for the set $\{1,...,n\}$.

\section{Markov Decision Process (MDP)}\label{sec:MDP}
\subsection{Setting and notations}
An MDP is defined by a tuple $(\X,\A,\bm{P},\bm{r},\bm{p}_{0},\lambda)$, where $\X$ is the set of states and $\A$ is the set of actions.  In this review we focus on finite sets of states and actions, and we write
 $|\X|=n < + \infty,  |\A|=A < + \infty.$ 
The kernel $\bm{P} \in \R^{n \times A \times n}$ models the transition rates from state-action pairs to the next states, and $\bm{r} \in \R^{n \times A}$ is the state-action reward.  There is an initial distribution $\bm{p}_{0} \in \R^{n}_{+}$ over the set of states,  and a discount factor $\lambda \in (0,1)$.

 A (stationary) \textit{policy} $\pi \in \left( \Delta(A) \right)^{n}$ maps each state to a probability distribution over the set of actions $\A$. For each policy $\pi$,  the \textit{value vector} $\bm{v}^{\pi} \in \R^{n}$  is defined as \[v^{\pi}_{s} = \E^{\pi, \boldsymbol{P}} \left[ \sum_{t=0}^{\infty} \lambda^{t}r_{s_{t}a_{t}} \; \bigg| \; s_{0} = s \right], \forall \; s \in \X,\] where $(s_{t},a_{t})$ is the state-action pair visited at time $t$. 
 From the dynamic programming principle, $\bm{v}^{\pi}$ satisfies the following recursion:
 \[ v^{\pi}_{s} =  \sum_{a \in \A} \pi_{sa}\left(r_{sa}+\lambda \sum_{s' \in \X} P_{sas'}v^{\pi}_{s'} \right) = \sum_{a \in \A} \pi_{sa}\left(r_{sa}+\lambda \bm{P}_{sa}^{\top}\bm{v}^{\pi} \right), \forall \; s \in \X,\]
 where $\bm{P}_{sa} \in \Delta(n)$ is the probability distribution over the next state given a state-action pair $(s,a)$.
We simply reformulate the previous equality as
 \[ \bm{v}^{\pi} = \left( \bm{I} - \lambda \bm{P}_{\pi} \right)^{-1}\bm{r}_{\pi},\]
 where $\bm{P}_{\pi} \in \R^{n \times n},\bm{r}_{\pi} \in \R^{n}$ are the transition matrix and the one-step expected reward vector induced by the policy $\pi$:
 \[ P_{\pi,ss'} = \sum_{a \in \A} \pi_{sa}P_{sas'}, r_{\pi,s} = \sum_{a \in \A} \pi_{sa}r_{sa}, \forall \; (s,s') \in \X \times \X.\]
 The goal of the decision-maker is to compute a policy $\pi^{*}$ that maximizes the infinite horizon expected discounted reward, defined as $R(\pi) = \bm{p}_{0}^{\top} \bm{v}^{\pi}.$ In this review we focus on Value Iteration (VI) and Policy Iteration (PI) algorithms to compute an optimal policy; we refer the reader to \cite{Puterman} for a detailed discussion about other algorithms.

\subsection{Value Iteration}\label{sec:VI}
Value Iteration was introduced by Bellman \citep{bellman1966dynamic}.
Define the \textit{Bellman operator} $T: \R^{n} \rightarrow \R^{n}$, where for $\bm{v} \in \R^{n}$,
\begin{equation}\label{eq:T_max-def}
T(\bm{v})_{s} = \max_{a \in \A} \{ r_{sa} + \lambda \cdot \bm{P}_{sa}^{\top}\bm{v} \}, \forall \; s \in \X.
\end{equation}
The operator $T$ is an order-preserving contraction of $\left(\R^{n},\| \cdot \|_{\infty} \right)$, where for $\bm{v},\bm{w} \in \R^{n}$, we have $ \bm{v} \leq \bm{w} \Rightarrow T(\bm{v})  \leq T(\bm{w}),$ and $\| T(\bm{v}) - T(\bm{w}) \|_{\infty}  \leq \lambda \cdot \| \bm{v} - \bm{w}\|_{\infty}.$ Crucially, note that $T$ is not differentiable everywhere, because of the $\max$ in its definition. The value iteration (VI) algorithm is defined as follows:
\begin{equation}\label{alg:VI}\tag{VI}
\bm{v}_{0} \in \R^{n}, \bm{v}_{t+1} = T(\bm{v}_{t}), \forall \; t \geq 0.
\end{equation}
The following theorem gives the convergence rate and stopping criterion for \ref{alg:VI}.
\begin{theorem}[\cite{Puterman}, Chapter 6.3]\label{th:conv-rate-VI-puterman}
\begin{enumerate}
\item The value vector $\bm{v}^{*}$ of the optimal policy $\pi^{*}$ is the unique fixed point of the operator $T$.
\item Let $\left(\bm{v}_{t}\right)_{t \geq 0}$ be generated by \ref{alg:VI}. For any $t \geq 0,$ we have $ \| \bm{v}_{t} - \bm{v}^{*} \|_{\infty} \leq \lambda^{t} \| \bm{v}_{0} - \bm{v}^{*} \|_{\infty}$. 
\item If
%\begin{equation}\label{eq:eps-approx-policy}
$\| \bm{v}_{t} - \bm{v}_{t+1} \|_{\infty} \leq \epsilon (1-\lambda) (2 \lambda)^{-1} $
%\end{equation}
then $ \| \bm{v}^{\pi_{t}} - \bm{v}^{*} \|_{\infty} \leq \epsilon,$
where $\bm{v}^{\pi_{t}}$ is the value vector of $\pi_{t}$, the policy attaining the maximum in each component of $T(\bm{v}_{t}).$
\end{enumerate}
\end{theorem}
Therefore, \ref{alg:VI} converges to $\bm{v}^{*}$ at a linear rate of $\lambda \in (0,1)$.
 However,  in many applications, the future rewards are almost as important as the reward in the current period and $\lambda$ is close to $1$.  In this case, the convergence of \ref{alg:VI} becomes very slow.
\paragraph{Value Computation.}
The problem of computing the value vector of a policy is also crucial, see for instance the \textit{Policy Evaluation} step in the Policy Iteration algorithm, presented in the next section.
The operator $T_{\pi}:\R^{n} \rightarrow \R^{n}$ associated with a policy $\pi$ is defined as
\begin{equation}\label{eq:T_pi}
T_{\pi}(\bm{v})_{s} = \sum_{a \in \A} \pi_{sa}\left(r_{sa} + \lambda \bm{P}_{sa}^{\top}\bm{v}\right), \forall \; s \in \X.
\end{equation}
Note that $T_{\pi}$ is an affine operator and a contraction for $\| \cdot \|_{\infty}$. The unique fixed point of $T_{\pi}$ is $\bm{v}^{\pi}$, the value vector of the policy $\pi$.
Therefore, the following algorithm is called \textit{Value Computation} (VC): 
\begin{equation}\label{alg:VC}\tag{VC}
\bm{v}_{0} \in \R^{n}, \bm{v}_{t+1} = T_{\pi}(\bm{v}_{t}), \forall \; t \geq 0.
\end{equation}
For the same reason as Algorithm \ref{alg:VI},  the sequence of vectors $\left(T_{\pi}^{t}(\bm{v}_{0})\right)_{t \geq 0}$ generated by \ref{alg:VC} converges linearly to $\bm{v}^{\pi}$ with a rate of $\lambda$, for any initial vector $\bm{v}_{0}$~\citep{Puterman}.  

\subsection{Policy Iteration}
Policy Iteration was developed by   Howard~\citep{howard-1960} and Bellman~\citep{bellman1966dynamic}.  The algorithm runs as follow. 
\begin{algorithm}[H]
\caption{Policy Iteration (PI)}\label{alg:PI}
\begin{algorithmic}[1]
\State \textbf{Initialize} $\pi_{0}$ at random.
\For{$t \geq 0$}
\State (\textit{Policy Evaluation)} Choose $ \bm{v}_{t}= \bm{v}^{\pi_{t}}$  the value vector of $\pi_{t}$.
\State (\textit{Policy Improvement) } Choose $\pi_{t+1}$ with $\pi_{t+1,s} \in \Delta(A)$ solving $\max_{\pi} \sum_{a \in \A} \pi_{a} \left( r_{sa} + \lambda \bm{P}_{sa}^{\top}\bm{v}_{t}\right), \forall \; s \in \X.$
\State (\textit{Stopping Criterion)} Stop when $\pi_{t+1}=\pi_{t}$.
\State Increment $t \leftarrow t+1.$
\EndFor
\end{algorithmic}
\end{algorithm}
In the Policy Improvement step,  Policy Iteration computes $\pi_{t+1}$ as a greedy one-step update,  given the value vector $\bm{v}^{\pi_{t}}$ related to $\pi_{t}$.
Prior to \cite{Ye}, the proofs of convergence of Policy Iteration \citep{howard-1960,bellman1966dynamic} relied on the fact that the policies $\pi_{1},...,\pi_{t}$ monotonically improve the objective function.  In particular,  the \textit{Policy Improvement} step guarantees that $\pi_{t+1}$ is a strictly better policy than $\pi_{t}$, unless $\pi_{t+1}=\pi_{t}$.  As there are only finitely many policies (in the case of finite numbers of states and actions sets), the Policy Iteration algorithm terminates. As a side note,  Policy Iteration is equivalent to a block-wise update rule when running the simplex algorithm on the linear programming reformulation of MDPs. 
 
 In Theorem 4.2 in \cite{Ye}, the author proves that the number of iterations of Policy Iteration is bounded from above by \[\dfrac{n^{2}A}{1-\lambda}\log\left(\dfrac{n^{2}}{1-\lambda}\right).\]
 Note that this bound does not depend of the order of magnitude of the rewards $r_{sa}$, so that the worst-case complexity of Policy Iteration is actually strongly polynomial (for fixed discount factor $\lambda$).  It is remarkable that even though the proofs presented in \cite{Ye} rely on simple tools from the theory of linear programming.
 These results have been further improved in \cite{scherrer2016improved} for MDPs and extended to two-player games in \cite{hansen2013strategy}.
\section{Value Iteration as a first-order method}\label{sec:VI-GD}
From Section \ref{sec:VI}, finding an optimal policy is equivalent to solving $\bm{v}=T(\bm{v})$.  It has been noted since as far back as \cite{bertsekas1995dynamic} that the operator $\bm{v} \mapsto \left(\bm{I} - T\right)(\bm{v})$ can be treated as the gradient of an unknown function.  We review here the results from \cite{goyal2019first}, where the authors build upon this analogy to define first-order methods for MDPs,  extend Nesterov's acceleration and Polyak's momentum to MDPs and present novel lower bounds on the performances of value iteration algorithms. We also review novel connections between Mirror Descent \citep{nemirovsky1983problem},  Primal-Dual Algorithm \citep{ChambollePock16}, and Value Iteration at the end of this section.  A review of the classical results for first-order methods in convex optimization can be found in Appendix \ref{sec:first-order}.
\subsection{First-order methods for MDPs}
Considering that $\bm{v} \mapsto \left(\bm{I} - T\right)(\bm{v})$ as the gradient of an unknown function, \cite{goyal2019first} define first-order methods for MDPs as follows.
\begin{definition} An algorithm is a \textit{first-order method} for MDPs if it produces a sequence of iterates $\left( \bm{v}_{t} \right)_{t \geq 0}$ satisfying \[ \bm{v}_{0}=\bm{0}, \bm{v}_{t+1} \in span \{\bm{v}_{0},...,\bm{v}_{t}, T(\bm{v}_{0}), ... T(\bm{v}_{t}) \}, \forall \; t \geq 0.\]
\end{definition}
Clearly, \ref{alg:VI} is a first-order method for MDP. However,  Policy Iteration is not, since it relies on the update $\bm{v}_{t} = \bm{v}^{\pi_{t}}$ (Policy Evaluation step).
\paragraph{Choice of parameters.}
Recall that for a differentiable, $\mu$-strongly convex, $L$-Lipschitz continuous function $f: \R^{n} \rightarrow \R$, the following inequalities hold: for all vectors $\bm{v},\bm{w} \in \R^{n},$
\begin{align*}
\mu \cdot \| \bm{v} - \bm{w} \|_{2} & \leq \| \nabla f(\bm{v}) - \nabla f (\bm{w}) \|_{2},  \\ 
\| \nabla f(\bm{v}) - \nabla f (\bm{w}) \|_{2} & \leq L \cdot \| \bm{v} - \bm{w} \|_{2}.
\end{align*}
In the case of the Bellman operator $T$,  for any vectors $\bm{v},\bm{w} \in \R^{n},$ the triangle inequality gives
\begin{align*}
(1-\lambda) \cdot \| \bm{v} - \bm{w} \|_{\infty} & \leq \| (\bm{I} - T)(\bm{v}) - (\bm{I} - T)(\bm{w}) \|_{\infty},
\\
\| (\bm{I} - T)(\bm{v}) - (\bm{I} - T)(\bm{w}) \|_{\infty} & \leq (1+\lambda) \cdot \| \bm{v} - \bm{w} \|_{\infty}. 
\end{align*}
 In convex optimization, the constant $\mu$ and $L$ can be used to tune the step sizes of the algorithms.  Because of the analogy between the last four equations,  the choice of $\mu = 1-\lambda$ and $L = 1+ \lambda$ will be used to tune the step sizes of the novel algorithms for MDPs, inspired from first-order methods in convex optimization. We also use the notation $\kappa = \mu/L$. 
\paragraph{Relaxed Value Iteration.}
For $\alpha > 0$, \cite{goyal2019first} consider the following \textit{Relaxed Value Iteration} (RVI) algorithm:
\begin{equation}\label{alg:R-VI}\tag{RVI}
\bm{v}_{0} \in \R^{n}, \bm{v}_{t+1} = \bm{v}_{t}-\alpha \left( \bm{v}_{t}- T(\bm{v}_{t}) \right), \forall \; t \geq 0.
\end{equation}
Note the analogy with the Gradient Descent algorithm \ref{alg:GD} (presented in Appendix \ref{sec:first-order}), when $ \bm{v} \mapsto (\bm{I}-T)(\bm{v})$ is the gradient of an unknown function.
\ref{alg:R-VI} is also considered in \cite{relax-VI-2,relax-VI-1} without an explicit connection to \ref{alg:GD}.   
 Formal convergence guarantees are provided in \cite{goyal2019first}; the main difficulty is the use of $\| \cdot \|_{\infty}$ instead of $\| \cdot \|_{2}$ in the properties that characterize the Bellman operator $T$. Recall that a sequence $\left( \bm{v}_{t} \right)_{t \geq 0}$ converges linearly to $\bm{v}^{*} \in \R^{n}$ at a rate of $\rho \in (0,1)$ if $\bm{v}_{t} = \bm{v}^{*} + o \left( \rho^{t} \right).$
\begin{proposition}[\cite{goyal2019first}, Proposition 1]\label{prop:VI-step size} Consider $\left( \bm{v}_{t} \right)_{t \geq 0}$ the iterates of \ref{alg:R-VI} and $\bm{v}^{*}$ the unique fixed point of the Bellman operator $T$.
Recall that $ \mu = 1-\lambda,  L=1+\lambda$ and $\kappa = \mu/L$.
\begin{enumerate}
\item The sequence $(\bm{v}_{t})_{t \geq 0}$ produced by \ref{alg:R-VI} does converge at a linear rate to $\bm{v}^{*}$ as long as $\alpha \in (0,2/L)$.
\item For $\alpha =1=2/(L+\mu)$, \ref{alg:R-VI} converges linearly to $\bm{v}^{*}$ at a rate of $\lambda=(1-\kappa)/(1+\kappa)$.
\end{enumerate}
\end{proposition} Note the analogy with the classical convergence results for \ref{alg:GD} (as presented in Proposition \ref{prop:GD} in Appendix \ref{sec:first-order}).  In particular, \ref{alg:VI} is a first-order method for MDP which converges to $\bm{v}^{*}$ at a rate of $(1-\kappa)/(1+\kappa)$, where $\kappa = \mu/L = (1-\lambda)/(1+\lambda).$ When $\lambda \rightarrow 1$, $\kappa \rightarrow 0$ and the convergence of \ref{alg:VI} becomes very slow. Also, the same results hold for \ref{alg:VC}, i.e., for a fixed policy $\pi$, \ref{alg:VC} converges to $\bm{v}^{\pi}$ at a rate of $(1-\kappa)/(1+\kappa)$.
\subsection{Nesterov's acceleration and Polyak's momentum for Value Iteration}
\paragraph{Nesterov Accelerated Value Iteration.}
Still considering $\bm{v}-T(\bm{v})$ as the gradient of some unknown function, it is possible to write an \textit{Accelerated Value Iteration} (AVI), building upon Accelerated Gradient Descent \ref{alg:AGD} (see Appendix \ref{sec:first-order} for the exact definition and convergence rates of \ref{alg:AGD}). In particular, \cite{goyal2019first} consider the following algorithm.
\begin{equation}\label{alg:AVI}\tag{AVI}
\bm{v}_{0},\bm{v}_{1} \in \R^{n}, \begin{cases}
    \bm{h}_{t}=\bm{v}_{t}+\gamma\cdot \left( \bm{v}_{t}-\bm{v}_{t-1} \right), \\
	\bm{v}_{t+1} \gets \bm{h}_{t}-\alpha \left( \bm{h}_{t}- T \left( \bm{h}_{t} \right) \right), \end{cases}\forall \; t \geq 1.
\end{equation}
The choice of step sizes for \ref{alg:AVI} is the same as for \ref{alg:AGD} in convex optimization (e.g., Theorem \ref{th:AGD-MGD} in Appendix \ref{sec:first-order}) but here $\mu= (1- \lambda), L=(1+ \lambda)$. This yields
\begin{equation}\label{eq:tuning-AVI}
 \alpha = 1/(1+\lambda),
\gamma = \left(1-\sqrt{1- \lambda^{2}}\right)/\lambda.
 \end{equation}
\paragraph{Polyak's Momentum Value Iteration.}
It is also possible to write a \textit{Momentum Value Iteration} algorithm (MVI):
\begin{equation}\label{alg:MVI}\tag{MVI}
\bm{v}_{0},\bm{v}_{1} \in \R^{n},
	\bm{v}_{t+1} =\bm{v}_{t}-\alpha \left( \bm{v}_{t}-T \left( \bm{v}_{t} \right)\right)  + \beta \cdot \left( \bm{v}_{t}-\bm{v}_{t-1}  \right) , \forall \; t \geq 1.
\end{equation}
The same choice as for Momentum Gradient Descent (\ref{alg:MGD}, see Theorem \ref{th:AGD-MGD} in Appendix \ref{sec:first-order}) gives the following step sizes:
 \begin{equation}\label{eq:tuning-MVI}
\alpha  =2/(1+\sqrt{1-\lambda^{2}}), \beta  =(1-\sqrt{1-\lambda^{2}})/(1+\sqrt{1-\lambda^{2}}).
\end{equation} 
\paragraph{Properties for non-affine operators.}
Unfortunately,  algorithms \ref{alg:AVI} and \ref{alg:MVI} may diverge on some MDP instances. In all generality, the analysis of \ref{alg:AVI} and \ref{alg:MVI} is related to the computation of the \textit{joint spectral radius} of a specific set of matrices, which is usually a hard problem \citep{blondel-NP-hard}.  \cite{vieillard2020momentum} obtain convergence of a momentum variant of $Q$-learning, under a strong assumption on the sequence of iterates (Assumption 1 in \cite{vieillard2020momentum}, closely related to the joint spectral radius of the visited matrices), and study this assumption empirically.   Despite this, the analysis of the case where $T$ is an affine operator reveals some interesting connections on the convergence properties of \ref{alg:AVI} and \ref{alg:MVI} compared to \ref{alg:AGD} and \ref{alg:MGD}.
\paragraph{Convergence rate for affine operators.}
Let $\pi$ be a policy and consider Algorithm \ref{alg:AVI} and Algorithm \ref{alg:MVI} where the Bellman operator $T$ is replaced with the Bellman recursion operator $T_{\pi}$, defined in \eqref{eq:T_pi}.  \cite{goyal2019first} refer to these new algorithms as \textit{Accelerated Value Computation} (AVC) and \textit{Momentum Value Computation} (MVC), in reference to Algorithm \ref{alg:VC}. Under conditions on the spectral radius of $\bm{P}_{\pi}$,  \cite{goyal2019first} provide the following convergence rates.
\begin{theorem}[\cite{goyal2019first}, Theorem 1 and Theorem 2]\label{th:AVC-MVC}
Assume that $\pi$ is a policy such that the spectrum of the transition matrix $\bm{P}_{\pi}$ is  contained in $\R$.
\begin{enumerate}
\item  Let $\alpha, \gamma$ be as in \eqref{eq:tuning-AVI}.  Then Algorithm AVC converges to $\bm{v}^{*}$ at a linear rate of $1- \sqrt{\kappa}$.
\item  Let $\alpha, \gamma$ be as in \eqref{eq:tuning-MVI}. 
 Then Algorithm MVC converges to $\bm{v}^{*}$ at a linear rate of $(1-\sqrt{\kappa})/(1+\sqrt{\kappa})$.
\end{enumerate}
\end{theorem}
Note the analogy between the convergence results for  Accelerated Gradient Descent and Momentum Gradient Descent (as presented in Theorem \ref{th:AGD-MGD} in Appendix \ref{sec:first-order}) and Theorem \ref{th:AVC-MVC} (for Accelerated Value Computation and Momentum Value Computation).  In particular, when $\lambda \rightarrow 1$ then $\kappa \rightarrow 0$ and both AVC and MVC enjoy significantly better convergence rates than \ref{alg:VC}, which achieves a convergence rate of $(1-\kappa)/(1+\kappa)$. The proof essentially relies on reformulating the AVC and MVC updates as second-order recursions on the sequence of iterates and inspecting the eigenvalues of the matrix $\bm{B}_{\pi}$ defining the dynamics of the iterates.  Under the assumption that the eigenvalues of $\bm{P}_{\pi}$ are real, the matrix $\bm{B}_{\pi}$ has maximum eigenvalue of at most $1-\sqrt{\kappa}$ in the case of AVC, and at most $(1-\sqrt{\kappa})/(1+\sqrt{\kappa})$ in the case of MVC.  Without this assumption, this spectral radius of the matrix $\bm{B}_{\pi}$ may be more than $1$, so that $(\bm{v}_{t})_{t \geq 0}$ may diverge (see Section 5.1.2 in \cite{goyal2019first}).  The proof in \cite{goyal2019first} is close to the methods in Section 3 in \cite{flammarion2015averaging}, pertaining to the analysis of algorithms from convex optimization such as \ref{alg:AGD}, \ref{alg:MGD} and linear averaging as second-order recursions.

The condition that the spectrum of $\bm{P}_{\pi}$ is contained in $\R$ is satisfied, in particular, when $\bm{P}_{\pi}$ is the transition matrix of a reversible Markov chain.  \cite{akian2020multiply} provide more general conditions on the spectrum of $\bm{P}_{\pi}$ for AVC to converge, and extend AVC to $d$-th order acceleration, where the inertial step involves the last $d$ iterates and the convergence rate approaches $1-\sqrt[d]{\kappa}$.
\subsection{Lower bound on any first-order method for MDP}
While in convex optimization, \ref{alg:AGD} achieves the best worst-case convergence rate over the class of smooth, convex functions,  this is not the case for \ref{alg:AVI} for MDPs. In particular, \cite{goyal2019first} prove the following lower-bound on any first-order method for MDPs.
\begin{theorem}[\cite{goyal2019first}, Theorem 3] \label{th:hard-instance} There exists an MDP instance $(\X,\A,\bm{P},\bm{r},\bm{p}_{0},\lambda)$ such that for any sequence of iterates $(\bm{v}_{t})_{t \geq 0}$ generated by a first-order method for MDPs,
the following lower bound holds for any step $t \in \{1,...,n-1\}$: \[\| \bm{v}_{t} - \bm{v}^{*} \|_{\infty} \geq \Omega \left( \lambda^{t} \right).\]
\end{theorem}
The proof of Theorem \ref{th:hard-instance} is constructive and relies on a hard MDP instance (which is neither reversible nor ergodic).
Note the analogy with the lower-bound on the convergence rates of any first-order methods in convex optimization (e.g., as presented in Theorem \ref{th:hard-instance-AGD} in Appendix \ref{sec:first-order}),  but recall that it is \ref{alg:VI} that converges at a linear rate of $\lambda$. Therefore, for MDPs, it is \textit{Value Iteration} that attains the lower bound on the worst-case complexity of first-order methods. This is in stark contrast with the analogous result in convex optimization (as presented in Theorem \ref{th:hard-instance-AGD}), where \ref{alg:AGD} attains the worst-case convergence rate.  It may be possible to obtain sharper lower-bounds for the convergence rate of first-order methods for MDPs by restraining the MDP instance to be reversible or ergodic.
\subsection{Connections with mirror descent and other first-order methods}
We review here some other algorithms from first-order convex optimization that have been recently used to design novel Value Iteration algorithms. 
\paragraph{Mirror Descent.}
The authors in \cite{geist2019theory} show that adding a regularization term in the Bellman update \eqref{eq:T_max-def} yields an algorithm resembling Mirror Descent (MD, \citep{nemirovsky1983problem}). 
In particular,  let $\bm{c}_{s,\bm{v}} = \left(r_{sa} + \lambda \bm{P}_{sa}^{\top}\bm{v}\right)_{a \in \A} \in \R^{n}$ for $\bm{v} \in \R^{n}$. Then the Bellman operator $T$ as in \eqref{eq:T_max-def} can rewritten as
$ T(\bm{v})_{s} = \max_{\pi \in \Delta(A)} \langle \pi,  \bm{c}_{s,\bm{v}} \rangle.$
Let $D: \Delta(A) \times \Delta(A) \rightarrow \R$  be a given Bregman divergence (for instance,  $D(\bm{x},\bm{y}) = \| \bm{x}-\bm{y}\|^{2}_{2}$ or the Kullback-Leibler divergence).  \cite{geist2019theory} define an algorithm resembling Mirror Descent and Value Iteration \ref{alg:OMD-VI} as follows:
\begin{equation}\label{alg:OMD-VI}\tag{MD-VI}
\bm{v}_{0} \in \R^{n},  \pi_{0} \text{ random },
\begin{cases} \pi_{t+1,s} \in \arg \max_{\pi \in \Delta(A)} \langle \pi,  \bm{c}_{s,\bm{v}_{t}} \rangle + D(\pi,\pi_{t,s}),  \forall \; s \in \X,\\
v_{t+1,s} = \langle \pi_{t+1,s},  \bm{c}_{s,\bm{v}_{t}} \rangle,  \forall \; s \in \X
\end{cases} \forall \; t \geq 0.
\end{equation}
It is also possible to replace the update on the value vector by 
\begin{equation}\label{eq:OMD-variant}
v_{t+1,s} = \langle \pi_{t+1,s},  \bm{c}_{\bm{v}_{t}} \rangle + D(\pi,\pi_{t,s}), \forall \; s \in \X, \forall \; t \geq 0,
\end{equation}
i.e., to include the regularization term in the update on $\bm{v}_{t+1}$.
%The following proposition gives intuition on the convergence of \ref{alg:OMD-VI}.
%\begin{proposition}[\cite{geist2019theory}, Corollary 3]
%Let $\left(\pi_{k}\right)$ the sequence of policies visited by \ref{alg:OMD-VI} (or its variants using the update \eqref{eq:OMD-variant}). Then 
%\[  \lim_{t \rightarrow + \infty}  \|  \bm{v}^{*} - \dfrac{1}{t}\sum_{\tau=1}^{t}\bm{v}^{\pi_{\tau}} \|_{\infty} = 0.\]
%\end{proposition}
%Note that the convergence is on the average of value vectors of the visited policies compared to $\bm{v}^{*}$ (i.e., the \textit{regret}),  and attains a rate of $o(1/t)$, compared to the linear convergence of Value Iteration.  
\cite{geist2019theory} provide the convergence rate of \ref{alg:OMD-VI} and a detailed discussion on the connection of \ref{alg:OMD-VI} and other classical algorithms for reinforcement learning, e.g., Trust Region Policy Optimization (TRPO, see Section 5.1 in \cite{geist2019theory}).
\paragraph{Other approaches.}
(Stochastic) Mirror Descent can also be used to solve the linear programming formulation of MDPs in min-max form \citep{jin2020efficiently,gong2020duality}.  
In \cite{grand2020first,grand2020scalable} the authors adapt the Primal-Dual Algorithm of \cite{ChambollePock16} for solving robust MDPs and distributionally robust MDPs.
While these papers provide new algorithms for solving MDPs,  they do not highlight a novel connection between classical algorithms for MDPs and classical algorithms for convex optimization.
The interested reader can find a concise presentation of connections between Franck-Wolfe algorithm \citep{frank1956algorithm},  Mirror Descent, dual averaging \citep{nesterov2009primal} and algorithms for reinforcement learning in the review of  \cite{vieillard2019connections}.  
We note that the algorithms for MDPs in \cite{vieillard2019connections} are optimizing in the space of policies ($\pi \in \left( \Delta(A) \right)^{n}$), thereby relying on methods from \textit{constrained} convex optimization. In contrast, in this review we consider that Value Iteration and Policy Iteration are optimizing over the space of value vectors ($\bm{v} \in \R^{n}$),  thereby relying on methods from \textit{unconstrained} convex optimization.
\section{Policy Iteration as a second-order method}\label{sec:PI-newton}
\subsection{Newton-Raphson with the Bellman operator}
The relation between Policy Iteration and Newton-Raphson method dates as far back as \cite{kalaba1959nonlinear} and \cite{pollatschek1969algorithms}.  We adapt here the presentation of \cite{puterman1979convergence}.  A review of the classical results for second-order methods in convex optimization can be found in Appendix \ref{sec:second-order}.
\paragraph{Jacobian of the gradient operator.} We introduce the following notations, which greatly simplify the exposition of the results.  Let us write $F:\bm{v} \mapsto \left( \bm{I} - T \right)(\bm{v})$.  Note that in Section \ref{sec:VI-GD}, we interpreted $F$ as the gradient of an unknown function $f: \R^{n} \rightarrow \R$. In order to develop \textit{second-order} methods,  one uses second-order information, i.e., the Hessian of the unknown function $f: \R^{n} \rightarrow \R$, which is the \textit{Jacobian} of the gradient operator $F: \R^{n} \rightarrow \R^{n}$. However, the map $F$ may not be differentiable at any vector $\bm{v}$. This is because the Bellman operator itself is not necessarily differentiable, as the maximum of some linear forms.  Still, it is possible to give a closed-form expression for the Jacobian of $F$, where it is differentiable. In particular, we have the following lemma.
\begin{lemma}\label{lem:jacobian-of-F}
Recall that $F(\bm{v}) = \left( \bm{I} - T\right)(\bm{v})$ and assume that $F$ is differentiable at $\bm{v} \in \R^{n}$.  Let $\bm{J}_{\bm{v}}$ be the Jacobian of $F$ at $\bm{v}$: \[\bm{J}_{\bm{v}} \in \R^{n \times n}, J_{\bm{v},ij} = \partial_{j} F_{i}(\bm{v}), \forall \; (i,j) \in [n] \times [n].\]
With this notation, we have  $\bm{J}_{\bm{v}}= \bm{I}-\lambda\bm{P}_{\pi(\bm{v})} $, where $\pi(\bm{v})$ attains the $\max$ in $T(\bm{v})$. 
Additionally, $\bm{J}_{\bm{v}}$ is invertible.
\end{lemma}
The fact that the matrix $\bm{J}_{\bm{v}}$ is invertible directly follows from $\lambda \in (0,1)$ and $\bm{P}_{\pi(\bm{v})}$ being a stochastic matrix.
\paragraph{Analytic representation of Policy Iteration.} Let us now analyze the update from $\bm{v}_{t} = \bm{v}^{\pi_{t}}$ to $\bm{v}_{t+1} = \bm{v}^{\pi_{t+1}}$ in Policy Iteration.
In the policy improvement step,  $\pi_{t+1}$ is obtained as a row-wise solution of $\max_{\pi} \{ \bm{r}_{\pi} + \lambda \bm{P}_{\pi}\bm{v}^{\pi_{t}} \}.$  The vector $\bm{v}_{t+1}$ is then updated as the value vector $\bm{v}^{\pi_{t+1}}$ of $\pi_{t+1}$.
From the definition of the value vector,  we have the recursion 
\[\bm{v}^{\pi_{t+1}} = \bm{r}_{\pi_{t+1}} + \lambda \bm{P}_{\pi_{t+1}}\bm{v}^{\pi_{t+1}},\]
and we can invert this equation to obtain
$ \bm{v}^{\pi_{t+1}} = \left( \bm{I} - \lambda \bm{P}_{\pi_{t+1}} \right)^{-1} \bm{r}_{\pi_{t+1}}.$ 
Note that with the notation of Lemma \ref{lem:jacobian-of-F}, we have $\bm{v}^{\pi_{t+1}} = \bm{J}_{\bm{v}_{t}}^{-1}\bm{r}_{\pi_{t+1}}.$
Adding and subtracting $\bm{v}_{t}$ from the previous equality yields
%\begin{align}
% \bm{v}^{\pi_{t+1}} & = \left( \bm{I} - \lambda \bm{P}_{\pi_{t+1}} \right)^{-1} \bm{r}_{\pi_{t+1}} \nonumber \\
% & = \bm{v}_{t} - \bm{v}_{t} + \left( \bm{I} - \lambda \bm{P}_{\pi_{t+1}} \right)^{-1} \bm{r}_{\pi_{t+1}} \nonumber \\
% & = \bm{v}_{t} - \left( \bm{I} - \lambda \bm{P}_{\pi_{t+1}} \right)^{-1}\left( \bm{I} - \lambda \bm{P}_{\pi_{t+1}} \right)\bm{v}_{t} +  \left( \bm{I} - \lambda \bm{P}_{\pi_{t+1}} \right)^{-1} \bm{r}_{\pi_{t+1}} \nonumber \\
% & = \bm{v}_{t} - \left( \bm{I} - \lambda \bm{P}_{\pi_{t+1}} \right)^{-1} \left(-  \bm{r}_{\pi_{t+1}}+\left( \bm{I} - \lambda \bm{P}_{\pi_{t+1}} \right)\bm{v}_{t} \right) \nonumber \\
% & =  \bm{v}_{t} - \left( \bm{I} - \lambda \bm{P}_{\pi_{t+1}} \right)^{-1} \left(\bm{v}_{t} -   \bm{r}_{\pi_{t+1}} - \lambda \bm{P}_{\pi_{t+1}}\bm{v}_{t} \right) \nonumber \\
%  & =  \bm{v}_{t} - \left( \bm{I} - \lambda \bm{P}_{\pi_{t+1}} \right)^{-1} \left(\bm{v}_{t} -  T(\bm{v}_{t}) \right),\label{eq:newton-as-PI}
%\end{align}
\begin{align}
 \bm{v}^{\pi_{t+1}} & =\bm{J}_{\bm{v}_{t}}^{-1} \bm{r}_{\pi_{t+1}} \nonumber \\
 & = \bm{v}_{t} - \bm{v}_{t} + \bm{J}_{\bm{v}_{t}}^{-1} \bm{r}_{\pi_{t+1}} \nonumber \\
 & = \bm{v}_{t} -\bm{J}_{\bm{v}_{t}}^{-1}\bm{J}_{\bm{v}_{t}}\bm{v}_{t} +  \bm{J}_{\bm{v}_{t}}^{-1} \bm{r}_{\pi_{t+1}} \nonumber \\
 & = \bm{v}_{t} -\bm{J}_{\bm{v}_{t}}^{-1} \left(-  \bm{r}_{\pi_{t+1}}+\bm{J}_{\bm{v}_{t}}\bm{v}_{t} \right) \nonumber.
 \end{align}
 We can then replace $\bm{J}_{\bm{v}_{t}}$ by its expression from Lemma \ref{lem:jacobian-of-F} in the last part of the above equation to obtain
 \begin{align}
 \bm{v}^{\pi_{t+1}} & =  \bm{v}_{t} - \bm{J}_{\bm{v}_{t}}^{-1} \left( -   \bm{r}_{\pi_{t+1}} + \left( \bm{I} - \lambda \bm{P}_{\pi_{t+1}} \right) \bm{v}_{t} \right) \nonumber \\
 & =  \bm{v}_{t} - \bm{J}_{\bm{v}_{t}}^{-1} \left(\bm{v}_{t} -   \bm{r}_{\pi_{t+1}} - \lambda \bm{P}_{\pi_{t+1}}\bm{v}_{t} \right) \nonumber \\
  & =  \bm{v}_{t} -\bm{J}_{\bm{v}_{t}}^{-1} \left(\bm{v}_{t} -  T(\bm{v}_{t}) \right),\label{eq:newton-as-PI}
\end{align}
where \eqref{eq:newton-as-PI} follows from the definition of $\pi_{t+1}$ attaining the $\max$ in $T(\bm{v}_{t})$.
Equation \eqref{eq:newton-as-PI} brings down to the following theorem.
\begin{theorem}[ \cite{puterman1979convergence}, Theorem 2]\label{th:PI-as-newton}
Assume that $F$ is differentiable at any vector $\bm{v}_{t}$ visited by Policy Iteration.
With the notation of Lemma \ref{lem:jacobian-of-F}, the Policy Iteration algorithm can be written as
\[ \bm{v}_{t+1} = \bm{v}_{t} - \bm{J}_{\bm{v}_{t}}^{-1}F(\bm{v}_{t}), \forall \; t \geq 1.\]
\end{theorem}
This is exactly the update of the Newton-Raphson method in convex optimization for solving  $F(\bm{v})=\bm{0}$ (as presented in Appendix \ref{sec:second-order}).  Unfortunately,  the map $F$ may not be differentiable at \textit{any} $\bm{v}$, because of the maximization term in the definition of the Bellman operator $T$.  Therefore,  this reformulation does not immediately yields the convergence property and the convergence rate of Policy Iteration using the classical results from convex optimization (e.g., Theorem \ref{th:conv-rate-Newton} in Appendix \ref{sec:second-order}). \cite{puterman1979convergence} provide some assumptions on the Lipschitzness of $\bm{v} \mapsto \bm{J}_{\bm{v}}$ that lead to convergence of Policy Iteration, and this analysis was refined recently in \cite{cvetkovic2020greedy} for a variant of Policy Iteration. However, these assumptions are hard to verify on MDP instances; to the best of our knowledge, the most meaningful analysis of the convergence rate of Policy Iteration remains \cite{Ye} and its refinements \citep{hansen2013strategy,scherrer2016improved}.  Additional analysis of Policy Iteration based on value functions approximations (for Markov games) can be found in \cite{perolat2016softened}.
\subsection{Smoothing the Bellman operator}
Replacing the $\max$ in the Bellman operator with a log-sum-exp reformulation, \cite{kamanchi2019second} obtain convergence of a log-sum-exp version of Policy Iteration directly from the convergence of the Newton-Raphson method. In particular,  for $\beta>0$, define $T_{\beta} : \R^{n} \rightarrow \R^{n}$ such that
\begin{equation}\label{eq:softmax-Bellman}
T_{\beta}(\bm{v})_{s} = \dfrac{1}{\beta} \log \left( \sum_{a \in \A} \exp \left(\beta \left(r_{sa} + \lambda \bm{P}_{sa}^{\top}\bm{v}\right) \right)\right), \forall \; s \in \X.
\end{equation} 
The authors in \cite{kamanchi2019second} prove the following results. Note that the same type of results can be obtained by replacing the log-sum-exp transformation with a soft-max.
\begin{theorem}[\cite{kamanchi2019second}, Lemmas 2-3, Theorem 2]
\begin{enumerate}
\item The operator $T_{\beta}$ is a $\lambda$-contraction for $\| \cdot \|_{\infty}$.  In particular, $T_{\beta}$ has a unique fixed point $\bm{v}_{\beta}^{*} \in \R^{n}$.
\item When $\beta \rightarrow + \infty$, we have $\bm{v}_{\beta}^{*} \rightarrow \bm{v}^{*}$ where $\bm{v}^{*}$ is the fixed point of the Bellman operator $T$. In particular, $\| \bm{v}_{\beta}^{*} - \bm{v}^{*} \|_{\infty} \leq \beta^{-1} \cdot \lambda (1-\lambda)^{-1}\log(A).$
\item Let $\left(\bm{v}_{t}\right)_{t \geq 0}$ be the sequence of iterates generated by the Newton-Raphson method \ref{alg:newton-raphson} for solving the equation $\bm{v}-T_{\beta}(\bm{v})=\bm{0}$. Then $\left(\bm{v}_{t}\right)_{t \geq 0}$ converges to $\bm{v}_{\beta}^{*}$ at a quadratic rate, for any initial point $\bm{v}_{0}$.
\end{enumerate}
\end{theorem}
\section{A Quasi-Newton method for MDP: Anderson acceleration}\label{sec:anderson-mdp}
Specific applications of Anderson accelerations to MDPs has been considered in \cite{geist2018anderson} and \cite{zhang2020globally}.  We give here some motivations, the relation with Anderson acceleration in convex optimization and review the convergence results of \cite{zhang2020globally}. A review of the classical results for quasi-Newton methods in convex optimization can be found in Appendix \ref{sec:quasi-newton}.
\paragraph{Anderson acceleration for MDPs.} We adapt here the presentation of \cite{geist2018anderson}.
Fix an integer $m \in \N$.
In order to compute the next iterates $\bm{v}_{t+1}$, Anderson VI computes weights $\alpha_{t,0}, ..., \alpha_{t,m}$ and updates $\bm{v}_{t+1}$ as a linear combination of the last $(m+1)$-iterates $T(\bm{v}_{t}), ..., T(\bm{v}_{t-m})$:
\begin{equation}\label{eq:anderson-update}
\bm{v}_{t+1} = \sum_{i=0}^{m} \alpha_{t,i} T(\bm{v}_{t-m+i}).
\end{equation}

Crucially, the weights are updated to minimize the weighted residuals of the previous of the previous iterates, where the residual of a vector $\bm{v}$ is $F(\bm{v})=\bm{v}-T(\bm{v})$.  In particular, the vector $\bm{\alpha}_{t}$ is chosen to minimize the following least-square problem:
\begin{equation}\label{eq:optim-pb-update-alpha}
\min_{\bm{\alpha} \in \R^{m+1}, \bm{\alpha}^{\top}\bm{e}=1} \| \sum_{i=0}^{m} \alpha_{i} \left( \bm{v}_{t-m+i}-T(\bm{v}_{t-m+i}) \right) \|^{2}_{2}. 
\end{equation}
The program \eqref{eq:optim-pb-update-alpha} is better understood if $T$ is a linear operator (or close to being linear). In particular,  if $ \sum_{i=0}^{m} \alpha_{i} T(\bm{v}_{t-m+i}) \approx T\left(\sum_{i=0}^{m} \alpha_{i} \bm{v}_{t-m+i} \right)$,  then
\[  \| \sum_{i=0}^{m} \alpha_{i} \left( \bm{v}_{t-m+i}-T(\bm{v}_{t-m+i}) \right) \|_{2} \approx \| T\left(\sum_{i=0}^{m} \alpha_{i} \bm{v}_{t-m+i} \right) - \sum_{i=0}^{m} \alpha_{i} \bm{v}_{t-m+i}\|_{2}.\]
Therefore,  in this case the weights $\bm{\alpha}$ are chosen so that the next vector $\bm{v}(\bm{\alpha})=\sum_{i=0}^{m} \alpha_{i}  \bm{v}_{t-m+i}$ has a minimal $F(\bm{v}(\bm{\alpha}))$ in terms of $\ell_{2}$-norm, where $F(\bm{v}) = \left( \bm{I} - T \right)(\bm{v})$. The vector $\bm{v}_{t+1}$ is then chosen as $T(\bm{v}(\bm{\alpha}))$.
\paragraph{Relation to Anderson Acceleration in Convex Optimization.} 
We present here a reformulation of Anderson Value Iteration which highlights the connection to the Anderson updates  in convex optimization (as presented in Appendix  \ref{sec:quasi-newton}). The connection between Anderson acceleration and quasi-Newton methods was first shown in \cite{eyert1996comparative} and expanded in \cite{fang2009two,walker2011anderson,zhang2020globally}. We follow the lines of \cite{walker2011anderson} and \cite{zhang2020globally} here. 
Note that there is a closed-form solution of \eqref{eq:optim-pb-update-alpha} at period $t$ using variable eliminations.  In particular, the optimization program \eqref{eq:optim-pb-update-alpha} can be rewritten as
\[ \min_{\bm{\beta}\in \R^{m}} \| F(\bm{v}_{t}) - \Delta \bm{ \FF}_{t}\bm{\beta} \|_{2}^{2},\]
where $\Delta \bm{\FF}_{t} \in \R^{n \times m}$ is defined as
\begin{align*}
\Delta \bm{\FF}_{t} & = \left( \Delta F_{t-m+1}, ... \Delta F_{t} \right) \in \R^{n \times m}, \\
\Delta F_{t'} & =  F(\bm{v}_{t'}) - F(\bm{v}_{t'-1}) \in \R^{n},  \forall \; t' \geq 0,
\end{align*}
and where  $\bm{\beta} \in \R^{m}$ relates to $\bm{\alpha} \in \R^{m+1}$ as $\alpha_{1} = \beta_{1}, \alpha_{m+1} = 1-\beta_{m}$ and $\alpha_{i} = \beta_{i} - \beta_{i-1}$ otherwise. If $\Delta\bm{ \FF}_{t}$ is non-singular,  the solution $\bm{\beta}_{t}$ to the least-square problem \eqref{eq:optim-pb-update-alpha} is simply 
\[ \bm{\beta}_{t} = \left( \Delta \bm{\FF}_{t}^{\top}\Delta \bm{\FF}_{t}\right)^{-1}\Delta \bm{\FF}_{t}^{\top}F(\bm{v}_{t}).\]
Using the relation between $\bm{\alpha}$ and $\bm{\beta}$ and plugging this expression into the Anderson update \eqref{eq:anderson-update} yields
\begin{equation}\label{alg:AndersonVI-bad}\tag{AndVI-II}
\begin{aligned}
\bm{G}_{t} & = \bm{I} + \left(\Delta \bm{\VV}_{t} - \Delta \bm{\FF}_{t} \right) \left(\Delta \bm{\FF}_{t}^{\top}\Delta \bm{\FF}_{t}\right)^{-1}\Delta \bm{\FF}_{t}^{\top},\\
\bm{v}_{t+1} & = \bm{v}_{t} - \bm{G}_{t}F(\bm{v}_{t}),
\end{aligned}
\end{equation}
with $\Delta \bm{\VV}_{t} = \left( \bm{v}_{t-m+1}-\bm{v}_{t-m}, ..., \bm{v}_{t}-\bm{v}_{t-1} \right) \in \R^{n \times m}.$ The update on the sequence $\left(\bm{G}_{t}\right)_{t \geq 0}$ is the same as the \textit{Anderson update of type-II} in convex optimization (see \eqref{eq:bad-anderson-update} in Appendix \ref{sec:quasi-newton}).  In particular, $\bm{G}_{t}$ achieves the minimum of $\bm{G}_{t} \mapsto \| \bm{G}_{t} - \bm{I} \|_{F}$, under the conditions that $\bm{G}_{t}\Delta \bm{\FF}_{t} = \Delta \bm{\VV}_{t}$; we write $\| \cdot \|_{F}$ for the Frobenius norm on matrices.  This last equation can be interpreted as an inverse multi-secant condition, since \[\bm{G}_{t}\Delta \bm{\FF}_{t} = \Delta \bm{\VV}_{t} \iff \bm{G}_{t} \left(F(\bm{v}_{t-m+i}) - F(\bm{v}_{t-m+i-1})\right) = \bm{v}_{t-m+i} -\bm{v}_{t-m+i-1}, \forall \; i=1,...,m.\]

The authors in \cite{zhang2020globally} also consider another version of Anderson Value Iteration,  using the  \textit{Anderson update of type-I} (see \eqref{eq:good-anderson-update} in Appendix \ref{sec:quasi-newton}):
\begin{equation}\label{alg:AndersonVI-good}\tag{AndVI-I}
\begin{aligned}\bm{J}_{t} & = \bm{I} + \left(\Delta \bm{\FF}_{t} - \Delta \bm{\VV}_{t} \right) \left(\Delta \bm{\VV}_{t}^{\top}\Delta \bm{\VV}_{t}\right)^{-1}\Delta \bm{\VV}_{t}^{\top}, \\
\bm{v}_{t+1} & = \bm{v}_{t} - \bm{J}_{t}^{-1}F(\bm{v}_{t}).
\end{aligned}
\end{equation}
The update \ref{alg:AndersonVI-good} approximates the Jacobian $\bm{J}_{\bm{v}_{t}}$ with $\bm{J}_{t}$, whereas \ref{alg:AndersonVI-bad} directly approximates $\bm{J}_{\bm{v}_{t}}^{-1}$ with $\bm{G}_{t}$.  In particular, $\bm{J}_{t}$ achieves the minimum of $\bm{J}_{t} \mapsto \| \bm{J}_{t} - \bm{I} \|_{F}$, under the condition that $\bm{J}_{t} \Delta \bm{\VV}_{t} = \Delta \bm{\FF}_{t}$.  It is possible to show that $\bm{G}_{t}$ and $\bm{J}_{t}$ are \textit{rank-one} updates from $\bm{G}_{t-1}$ and $\bm{J}_{t-1}$ \citep{zhang2020globally}.
The analysis of Anderson Value Iteration as in Algorithm \ref{alg:AndersonVI-good} and Algorithm \ref{alg:AndersonVI-bad} is complex.  In particular,  the  operator $T$ is not differentiable, so that it is not straightforward that \ref{alg:AndersonVI-good} or \ref{alg:AndersonVI-bad} converge.  
In order to prove convergence of Algorithm \ref{alg:AndersonVI-good},  the authors in \cite{zhang2020globally}:
\begin{itemize}
\item tackle the potential singularity of the matrices $\left( \Delta \bm{\VV}_{t}^{\top}\Delta \bm{\VV}_{t} \right)_{t \geq 0}$ (with \textit{Powell-type} regularization),
\item tackle the potential singularity of the matrices $\left(\bm{J}_{t}\right)_{t \geq 0}$ (with \textit{restart checking}),
\item ensure sufficient decrease in $\| T(\bm{v}_{t}) - \bm{v}_{t} \|_{2}$ at every iteration (with \textit{safeguarding steps}).
\end{itemize} 
\cite{zhang2020globally} present a \textit{stabilized} version of the vanilla Anderson algorithm \ref{alg:AndersonVI-good}. This results in $\lim_{t \rightarrow + \infty} \bm{v}_{t} = \bm{v}^{*}$,  but the convergence rate is not known (Theorem 3.1 in \cite{zhang2020globally}), even though the algorithm enjoys good empirical performances, typically outperforming \ref{alg:VI}.  Note that this is the first result on the convergence of Anderson Acceleration,  \textit{without differentiability} of the operator $T$ (and with fixed-memory, i.e., with $m$ fixed).
There are many other quasi-Newton methods aside from Anderson acceleration, see Appendix \ref{sec:quasi-newton} for a short review of some classical methods. Note that in the practical implementation of Anderson acceleration, one typically chooses $1 \leq m \leq 5$. Therefore, algorithms based on information on only the last two iterates may still perform well.
In particular, it could be possible to develop algorithms for MDPs, based on Broyden or BFGS updates.  The convergence of \ref{alg:AndersonVI-bad} for non-differentiable operators (or a stabilized version of \ref{alg:AndersonVI-bad}) also remains an open question.
\small
\bibliographystyle{plainnat} 
\bibliography{FOM_RMDP}
\normalsize
\appendix
\section{A short review of gradient-based methods in convex optimization}\label{sec:convex-opt}
In this section we review the classical definitions, algorithms and convergence results from convex optimization (mostly referring to~\cite{boyd-2004,nesterov-book}). As the scope of gradient methods in convex optimization is extremely large, we only introduce here the notions that are relevant to our connections with algorithms for MDPs.

\subsection{First-order methods}\label{sec:first-order}

In this section we let $f: \R^{n} \rightarrow \R$ be a differentiable function.  We start with the following definitions of first-order algorithms for convex optimization and linear convergence.
\begin{definition}
A \textit{first-order algorithm} produces a sequence $\left( \bm{x}_{t} \right)_{t \geq 0}$ of points of $\R^{n}$ that only exploits first-order information at the previous iterates to compute the next iterate.  In particular, a first-order algorithm satisfies the following condition:
\[\bm{x}_{t+1} \in \bm{x}_{0} + span \{\nabla f (\bm{x}_{0}),...,\nabla f (\bm{x}_{t}) \}, \forall \; t \geq 0.\]
\end{definition}
We also use the following definition of linear convergence in convex optimization.
\begin{definition}
The convergence of $\left(\bm{x}_{t} \right)_{t \geq 0}$ to $\bm{x}_{f}$ a minimizer of $f$ is called \textit{linear} of rate $\rho \in (0,1)$ if \[\bm{x}_{t} = \bm{x}_{f} + o \left( \rho^{t} \right).\]
%$\| \bm{x}_{t} - \bm{x}_{f} \|_{2} \leq O \left( \rho^{t} \| \bm{x}_{0} - \bm{x}_{f} \|_{2} \right), \forall \; t \geq 1.$
\end{definition}

\paragraph{Gradient descent.}
The gradient descent algorithm (Algorithm \ref{alg:GD}) with fixed step size $\alpha$ is described as:
\begin{equation}\label{alg:GD}\tag{GD}
\bm{x}_{0} \in \R^{n}, \bm{x}_{t+1} = \bm{x}_{t}-\alpha \nabla f (\bm{x}_{t}), \forall \; t \geq 0.
\end{equation}
The intuition for \ref{alg:GD} comes from first-order Taylor expansions. In particular,  if $f$ is differentiable, we can write $f(\bm{x}+\bm{h}) = f(\bm{x}) + \bm{h}^{\top}\nabla f(\bm{x}) + o(\| \bm{h}\|_{2})$. Now using the update rule for $\bm{x}_{t+1}$ we obtain
\[ f(\bm{x}_{t+1}) = f(\bm{x}_{t}-\alpha \nabla f (\bm{x}_{t})) = f(\bm{x}_{t}) - \alpha \| \nabla f (\bm{x}_{t}) \|_{2}^{2} + o(\| \nabla f (\bm{x}_{t}) \|_{2}).\]
Therefore, for the right choice of step size $\alpha$ we can obtain $f(\bm{x}_{t+1}) \leq f(\bm{x}_{t})$. In particular, the following proposition gives the rate of convergence of \ref{alg:GD}. 
\begin{proposition}[\cite{nesterov-book}, Chapter 2.1.5)]\label{prop:GD}
Assume that $f$ is $\mu$-strongly convex and $L$-Lipschitz continuous. Let $\bm{x}^{*}$ be the minimizer of $f$ and $\kappa=\mu/L$.
\begin{enumerate}
\item The sequence $(\bm{x}_{t})_{t \geq 0}$ produced by \ref{alg:GD} does converge at a linear rate to $\bm{x}^{*}$ as long as $\alpha \in (0,2/L)$.
\item For $\alpha =2/(L+\mu)$,  \ref{alg:GD} converges linearly to $\bm{x}^{*}$ at a rate of $(1-\kappa)/(1+\kappa)$.
\item In both the previous cases, \ref{alg:GD} is a descent algorithm: $f(\bm{x}_{t+1}) \leq f(\bm{x}_{t}), \forall \; t \geq 1.$
\end{enumerate}
\end{proposition}
Note that they are many other step sizes strategy,  e.g.,  (backtracking) line search or choosing $\alpha_{t} = 1/t$ \citep{boyd-2004}. 
\paragraph{Acceleration and Momentum in convex optimization.}
Two popular first-order algorithms building upon \ref{alg:GD} are \textit{Accelerated Gradient Descent} (\ref{alg:AGD}) and \textit{Momentum} Gradient Descent (\ref{alg:MGD}). In particular, 
 Accelerated Gradient Descent \citep{nesterov-1983,nesterov-book}  builds upon \ref{alg:GD} by adding an intermediate extrapolation step:
\begin{equation}\label{alg:AGD}\tag{AGD}
\bm{x}_{0},\bm{x}_{1} \in \R^{n}, \begin{cases}
    \bm{h}_{t}=\bm{v}_{t}+\gamma\cdot \left( \bm{x}_{t}-\bm{x}_{t-1} \right), \\
	\bm{x}_{t+1} \gets \bm{h}_{t}-\alpha\nabla f (\bm{h}_{t}), \end{cases} \forall \; t \geq 1.
\end{equation}
In contrast, Momentum Gradient Descent builds upon \ref{alg:GD} by adding a momentum term in the update:
\begin{equation}\label{alg:MGD}\tag{MGD}
\bm{x}_{0},\bm{x}_{1} \in \R^{n},
	\bm{x}_{t+1} =\bm{x}_{t}-\alpha \nabla f (\bm{x}_{t})  + \beta\cdot \left( \bm{x}_{t}-\bm{x}_{t-1} \right) , \forall \; t \geq 1.
\end{equation}

The following theorem presents the rates of convergence of \ref{alg:AGD} and \ref{alg:MGD} for fixed step sizes. 

\begin{theorem}[\cite{heavyball-2015}, \cite{nesterov-book}]\label{th:AGD-MGD}
Assume that $f$ is $\mu$-strongly convex and $L$-Lipschitz continuous. Let $\bm{x}^{*}$ be the minimizer of $f$ and $\kappa=\mu/L$.
\begin{enumerate}
\item  Let $\alpha=1/L,\gamma = (\sqrt{L}-\sqrt{\mu})/(\sqrt{L}+\sqrt{\mu})$.  Then \ref{alg:AGD} converges to $\bm{x}^{*}$ at a linear rate of $1- \sqrt{\kappa}$.
\item Additionally,  assume that $f$ is twice continuously differentiable.  Let 
 $\alpha=4/(\sqrt{L} + \sqrt{\mu})^{2},\beta = (\sqrt{L}-\sqrt{\mu})^{2}/(\sqrt{L}+\sqrt{\mu})^{2}$. 
 Then \ref{alg:MGD} converges to $\bm{x}^{*}$ at a linear rate of $(1-\sqrt{\kappa})/(1+\sqrt{\kappa})$.
\end{enumerate}
\end{theorem}
When $\kappa$ is small, both \ref{alg:AGD} and \ref{alg:MGD} enjoy better convergence rates than $(1-\kappa)/(1+\kappa)$, the convergence rate of \ref{alg:GD}.
Note that the convergence of \ref{alg:MGD} relies on stronger assumption than the convergence of \ref{alg:AGD}. In particular, \ref{alg:MGD} may diverge if $f$ is not twice continuously differentiable \citep{heavyball-2015}. Additionally, \ref{alg:AGD} and \ref{alg:MGD} are \textit{not} descent algorithms: they do not necessarily produce estimates that result in a monotonically decreasing objective function. The objective value might increase for a few periods, before significantly decreasing afterward. This is known as \textit{oscillations}, and we refer to \cite{o2015adaptive} for a detailed study of the oscillation effects of \ref{alg:AGD}.
\paragraph{Lower bounds on the convergence rates.}

 \cite{nesterov-book} provides lower-bounds on the convergence rate of any first-order algorithm on the class of smooth, convex functions and on the class of smooth, strongly-convex functions. In particular, we recall the results for lower-bounds of first-order algorithms in smooth convex optimization.
 \cite{nesterov-book} considers $n = + \infty$, i.e., considers the space $\R^{\N}$ of infinite sequences of scalars, and proves the following lower-bound for smooth, strongly-convex functions.
\begin{theorem}[ \cite{nesterov-book}, Th. 2.1.13]\label{th:hard-instance-AGD}
For any $\bm{x}_{0} \in \R^{\N}$ and any step $t \geq 0$, there exists a $\mu$-strongly convex, $L$-Lipschitz continuous function $f: \R^{\N} \rightarrow \R$ such that for any sequence of iterates generated by a first-order method, the following lower bound holds :
\begin{align*}
\| \bm{x}_{t} - \bm{x}^{*} \|_{2}  \geq \left(\dfrac{1-\sqrt{\kappa}}{1+\sqrt{\kappa}} \right)^{t} \cdot \| \bm{x}_{0} - \bm{x}^{*} \|_{2}.
\end{align*}
\end{theorem}
The proof of Theorem \ref{th:hard-instance-AGD} relies on designing a hard instance. Additionally, \cite{nesterov-book} proves that the rate of convergence of AGD matches this lower-bound: \ref{alg:AGD} achieves the optimal rate of convergence over the class of smooth, convex functions, as well as over the class of smooth, strongly-convex functions.  Note that \ref{alg:MGD} achieves a better worst-case convergence rate than \ref{alg:AGD} but requires that $f$ is \textit{twice} continuously differentiable. In particular, \ref{alg:MGD} may diverge if $f$ is only continuously differentiable once \citep{heavyball-2015}.
\subsection{Second-order methods}\label{sec:second-order}
We focus here on the Newton-Raphson method (NR). In particular, consider solving the equation $F(\bm{x})=\bm{0}$, for a function $F: \R^{n} \rightarrow \R^{n}$.  If the function $F$ is differentiable, the \textit{Jacobian} $\bm{J}_{\bm{x}}$ of $F$ is defined for any $\bm{x} \in \R^{n}$, where $\bm{J}_{\bm{x}} \in \R^{n \times n}, J_{\bm{x},ij} = \partial_{j} F_{i}(\bm{x}), \forall \; (i,j) \in [n] \times [n]$.  The NR algorithm is:
\begin{equation}\label{alg:newton-raphson}\tag{NR}
\bm{x}_{0} \in \R^{n}, \bm{x}_{t+1} = \bm{x}_{t} -   \bm{J}_{\bm{x}_{t}}^{-1} F (\bm{x}_{t}), \forall \; t \geq 0.
\end{equation}
The intuition is as follows. If we linearize the condition $F(\bm{x})=\bm{0}$ close to $\bm{x}$, we have
 $F(\bm{x}+\bm{v}) \approx F(\bm{x}) + \bm{J}_{\bm{x}}\bm{v}$, 
 and $F(\bm{x}) + \bm{J}_{\bm{x}}\bm{v} = \bm{0} \iff \bm{v} = - \bm{J}_{\bm{x}}^{-1} F (\bm{x} ).$ Therefore,  for $\bm{x}$ close to $\bm{x}^{*}$,   the point $\bm{x} -  \bm{J}_{\bm{x}}^{-1} F (\bm{x})$ should be a good approximation of $\bm{x}^{*}$.
Under some conditions,  Algorithm \ref{alg:newton-raphson} converges to $\bm{x}^{*}$ a solution to $F(\bm{x})=\bm{0}$ and the convergence is quadratic.

\begin{theorem}\label{th:conv-rate-Newton}[\cite{boyd-2004} Section 9.5, \cite{poczos-lecture-convex-raphson}]
 Assume that:
 \begin{itemize}
 \item $\bm{J}_{\bm{x}} \succ \mu \bm{I}, \forall \; \bm{x} \in \R^{n}$. 
 \item $\exists \; M \geq 0, \| \bm{J}_{\bm{x}}- \bm{J}_{\bm{y}} \|_{2} \leq M \| \bm{x}-\bm{y} \|_{2}, \forall \; \bm{x},\bm{y} \in \R^{n}.$
 \end{itemize}
 Then there is a unique solution $\bm{x}^{*}$ to $F(\bm{x})=\bm{0}$.  Additionally,  if $\bm{x}_{0}$ is sufficiently close to $\bm{x}^{*}$: $\| \bm{x}_{0} - \bm{x}^{*} \|_{2} \leq 2\mu/3M$, then there exists $\rho > 0$, such that Algorithm \ref{alg:newton-raphson} enjoys a quadratic convergence rate to $\bm{x}^{*}$:
\[ \| \bm{x}_{t+1}-\bm{x}^{*} \|_{2} \leq \rho \| \bm{x}_{t} - \bm{x}^{*} \|_{2}^{2}.\]
\end{theorem}
In the case of quadratic convergence, the number of correct digits of $\bm{x}_{t}$ (compared to the optimal solution $\bm{x}^{*}$) roughly doubles at every iteration of the Newton-Raphson method (see Equation (9.35) in \cite{boyd-2004}). However, the convergence result of Theorem \ref{th:conv-rate-Newton} is very subtle. In particular,  if $\bm{x}_{0}$ is not close enough to $\bm{x}^{*}$,  \ref{alg:newton-raphson} may diverge. Additionally, the quadratic convergence result may hold only for a fairly small region around $\bm{x}^{*}$. Without the assumptions presented in Theorem \ref{th:conv-rate-Newton}, the algorithm may cycle or converge only in linear time.  We refer the reader to \cite{poczos-lecture-convex-raphson} for some examples of problematic cases. We also refer the reader to Chapter 7, 10 and 12 of \cite{ortega2000iterative} for a comprehensive analysis of Newton-Raphson method and its variations.
Note that it is possible to start in a \textit{damped phase} where we add a step size $\alpha_{t}$, i.e., to start with $\bm{x}_{t+1} = \bm{x}_{t} -   \alpha_{t} \bm{J}_{\bm{x}_{t}}^{-1} F (\bm{x}_{t})$ for adequate choices of the step sizes $\alpha_{t}$.  During the damped phase, \ref{alg:newton-raphson} has linear convergence, before transitioning to quadratic convergence when the iterates are close enough to the optimal solution $\bm{x}^{*}$, see \cite{boyd-2004}, Section 9.5.1.
\paragraph{Relation to Newton's method in convex optimization.}
We highlight here the relation between \ref{alg:newton-raphson} and Newton's method.
Let $f$ be a $\mu$-strongly convex and twice-differentiable function.  In this case, the Hessian $\nabla^{2} f (\bm{x})$ is defined at any point $\bm{x} \in \R^{n}$ and $\nabla^{2} f (\bm{x}) \succ \mu \bm{I}$. Newton's method to find a minimizer of $f$ follows from minimization the second-order Taylor expansion of $f$ at the current iterate and runs as follows:
 \begin{equation}\label{alg:netwon-convex}\tag{Newton}
\bm{x}_{0} \in \R^{n}, \bm{x}_{t+1} = \bm{x}_{t} -  \alpha_{t} \left(\nabla^{2} f (\bm{x}_{t}) \right)^{-1} \nabla f (\bm{x}_{t}), \forall \; t \geq 0. 
\end{equation}
The relation with \ref{alg:newton-raphson} comes from first-order optimality conditions.
In order to minimize $f$, we want to solve the system of equations $\nabla f (\bm{x}) = \bm{0}$.  We can then use \ref{alg:newton-raphson} on this equation in $\bm{x} \in \R^{n}$.
The operator $\bm{x} \mapsto \nabla f (\bm{x})$ is differentiable (since $f$ is twice-differentiable) and plays the role of the operator $F$ in the Newton-Raphson algorithm. The first assumption in Theorem \ref{th:conv-rate-Newton} is equivalent to $f$ being $\mu$-strongly convex and the second assumption in Theorem \ref{th:conv-rate-Newton} translates as $\exists \; M \geq 0, \| \nabla^{2}f(\bm{x}) - \nabla^{2}f(\bm{y}) \|_{2} \leq M \| \bm{x}-\bm{y} \|_{2}, \forall \; \bm{x},\bm{y} \in \R^{n}.$
\subsection{Quasi-Newton methods}\label{sec:quasi-newton}
Invented mostly in the 1960s and 1970s, quasi-Newton methods are still very popular to solve nonlinear optimization problems and nonlinear equations; their applications are very diverse and range from machine learning \citep{zhang2020globally} to quantum chemistry \citep{rohwedder2011analysis}. We describe here the \textit{Broyden updates} for constructing successive approximations of the Jacobians,  discuss their convergence rates and their relations to the BFGS update. We also detail \textit{Anderson acceleration}, as it relates to Anderson Value Iteration for MDPs (see Section \ref{sec:anderson-mdp}).

Let us consider again solving $F(\bm{v})=\bm{0}$ with $F: \R^{n} \rightarrow \R^{n}$.
Generally speaking, any method that replaces the exact computation of the (inverse of the) Jacobian matrices in the Newton-Raphson algorithm \ref{alg:newton-raphson} with an approximation is a \textit{quasi-Newton} method. 
This can be useful when the derivative of $F$ is hard to compute,  or the dimension of the problem is very large.  
In particular, a quasi-Newton method constructs a sequence of iterates $\left( \bm{x}_{t} \right)_{t \geq 0}$ and a sequence of matrices $\left(\bm{J}_{t}\right)_{t \geq 0}$ such that $\bm{J}_{t}$ is an approximation of the Jacobian $\bm{J}_{\bm{x}_{t}}$ for any $t \geq 0$ and
\begin{equation}\label{alg:quasi-newton}
\bm{x}_{0} \in \R^{n}, \bm{x}_{t+1} = \bm{x}_{t} -   \bm{J}_{t}^{-1} F (\bm{x}_{t}), \forall \; t \geq 0.
\end{equation}
\paragraph{Good and bad Broyden updates.}
Recall that the motivation for the Newton-Raphson algorithm \ref{alg:newton-raphson} comes from linearization:
\begin{equation}\label{eq:linearization-F}
F(\bm{x}_{t}+\bm{v}) \approx F(\bm{x}_{t}) + \bm{J}_{\bm{x}_{t}}\bm{v}.
\end{equation}
Therefore, when we compute the approximation $\bm{J}_{t}$ we want to maintain \eqref{eq:linearization-F}, which yields the so-called \textit{secant condition} on $\bm{J}_{t}$:
\begin{equation}\label{eq:secant-condition}
\bm{J}_{t}\Delta \bm{x}_{t}= \Delta F_{t}.
\end{equation}
for $\Delta \bm{x}_{t} = \bm{x}_{t} - \bm{x}_{t-1} , \Delta F_{t}  = F(\bm{x}_{t}) - F(\bm{x}_{t-1})$. Another common condition is that $\bm{J}_{t}$ and $\bm{J}_{t-1}$ coincides on the direction orthogonal to $\Delta \bm{x}_{t}$. This is the so-called \textit{no-change condition}:
\begin{equation}\label{eq:no-change-cond}
\bm{J}_{t}\bm{q} = \bm{J}_{t-1}\bm{q}, \forall \; \bm{q} \in \R^{n} \text{ such that } \bm{q}^{\top}\Delta \bm{x}_{t} = 0.
\end{equation}
There is a unique matrix satisfying \eqref{eq:secant-condition}-\eqref{eq:no-change-cond} and this yields the \textit{type-I} Broyden's update:
\begin{equation}\label{eq:broyden-good}\tag{Broyden-I}
\bm{J}_{t} = \bm{J}_{t-1} + \left( \Delta F_{t} - \bm{J}_{t-1} \Delta \bm{x}_{t} \right) \dfrac{\Delta \bm{x}_{t}^{\top}}{\Delta \bm{x}_{t}^{\top}\Delta \bm{x}_{t}}.
\end{equation}
Note that this is a rank one update from $\bm{J}_{t-1}$. Additionally,  \cite{dennis1977quasi} show that the update \eqref{eq:broyden-good}  minimizes $\bm{J}_{t} \mapsto \| \bm{J}_{t-1} - \bm{J}_{t} \|_{F}^{2}$ subject to the secant condition \eqref{eq:secant-condition}, where $\| \cdot \|_{F}$ is the Frobenius norm.
The update \eqref{eq:broyden-good} approximates the Jacobian and is sometimes called the \textit{good} Broyden's update, as it enjoys better empirical performances than the \textit{bad} Broyden update (\textit{type-II}),  which directly approximates the inverse of the Jacobian by considering both the secant condition \eqref{eq:secant-condition} and the no-change condition \eqref{eq:no-change-cond} in their inverse forms. In particular, the bad Broyden update finds a matrix $\bm{G}$ satisfying
\begin{equation}\label{eq:condition-bad-broyden}
\bm{G}_{t} \Delta F_{t} = \Delta \bm{x}_{t}, \bm{G}_{t}\bm{q} = \bm{G}_{t-1}\bm{q}, \forall \; \bm{q} \in \R^{n} \text{ such that } \bm{q}^{\top}\Delta F_{t}=0.
\end{equation}
This also has a closed-form solution:
\begin{equation}\label{eq:broyden-bad}\tag{Broyden-II}
\bm{G}_{t} = \bm{G}_{t-1} + \left( \Delta \bm{x}_{t} - \bm{G}_{t-1} \Delta F_{t} \right) \dfrac{\Delta \bm{x}_{t}^{\top} \bm{G}_{t-1}}{\Delta \bm{x}_{t}^{\top}\bm{G}_{t-1}\Delta \bm{x}_{t}},
\end{equation}
and $\bm{G}_{t}$ minimizes $\bm{G}_{t} \mapsto \| \bm{G}_{t-1} - \bm{G}_{t} \|_{F}^{2}$ subject to the secant condition in inverse form \eqref{eq:condition-bad-broyden}. Since $\bm{G}_{t}$ directly approximates the \textit{inverse} of the Jacobian, the update on $\bm{x}_{t}$ is $
\bm{x}_{t+1} = \bm{x}_{t} - \bm{G}_{t}F(\bm{x}_{t}).$ 
Note that there are many other possible approximation schemes for the Jacobians (or their inverses). We refer the reader to \cite{fang2009two}, Section 2 and Section 3, for more details.
\paragraph{Superlinear Convergence of Broyden's updates.}
Recall that the superlinear convergence of a sequence $\left(\bm{x}_{t} \right)_{t \geq 0}$ to $\bm{x}^{*}$ is equivalent to $\lim_{t \rightarrow \infty} \| \bm{x}_{t+1} - \bm{x}^{*} \|_{2} / \| \bm{x}_{t} - \bm{x}^{*} \|_{2} = 0$.
The result for the local superlinear convergence of the quasi-Newton iterations with Broyden's updates dates back to \cite{broyden1973local}. Crucially,  the next convergence results require  \textit{differentiable} operators.
\begin{theorem}[\cite{broyden1973local}, Theorems 3.2-3.4 ] \label{th:broyden-convergence}
Let $\bm{x}^{*} \in \R^{n}$ be such that $F(\bm{x}^{*}) = \bm{0}$.
Assume that:
\begin{itemize}
\item  $F: \R^{n} \rightarrow \R^{n}$ is differentiable on $\R^{n}$.
\item $\bm{J}_{\bm{x}^{*}}$ is non-singular.
\item $ \exists \; M \in \R,  \| \bm{J}_{\bm{x}} - \bm{J}_{\bm{x}^{*}} \|_{2} \leq M \| \bm{x} -\bm{x}^{*} \|_{2}, \forall \; \bm{x} \in \R^{n}.$
\end{itemize}
Then the following convergence results hold.
\begin{enumerate}
\item  Let $\left(\bm{x}_{t} \right)_{t \geq 0}$ following $\bm{x}_{t+1} = \bm{x}_{t} - \bm{J}_{t}^{-1}F(\bm{x}_{t})$ where $\left( \bm{J}_{t} \right)_{t \geq 0}$ is updated following \eqref{eq:broyden-good}.

If $\bm{x}_{0}$ is close enough to $\bm{x}^{*}$,  then $\lim_{t \rightarrow \infty} \bm{x}_{t} = \bm{x}^{*}$ and the convergence is superlinear.
\item  Let $\left(\bm{x}_{t} \right)_{t \geq 0}$ following $\bm{x}_{t+1} = \bm{x}_{t} - \bm{G}_{t}F(\bm{x}_{t})$ where $\left( \bm{G}_{t} \right)_{t \geq 0}$ is updated following \eqref{eq:broyden-bad}.

If $\bm{x}_{0}$ is close enough to $\bm{x}^{*}$,  then $\lim_{t \rightarrow \infty} \bm{x}_{t} = \bm{x}^{*}$ and the convergence is superlinear.
\end{enumerate}
\end{theorem}
\paragraph{Comparison to BFGS.}
Note that the updates \eqref{eq:broyden-good} and \eqref{eq:broyden-bad} are rank-one updates. In contrast, the \textit{Broyden-Fletcher-Goldfarb-Shanna} (BFGS) algorithm performs rank-two updates of the form
\begin{equation}\label{eq:BFGS}\tag{BFGS}
\bm{J}_{t} = \bm{J}_{t-1} + \dfrac{\Delta F_{t} \Delta F_{t}^{\top}}{\Delta F_{t}^{\top} \Delta \bm{x}_{t} }  - \dfrac{ \bm{J}_{t-1} \Delta \bm{x}_{t} \left( \bm{J}_{t-1} \Delta \bm{x}_{t} \right)^{\top}}{\Delta \bm{x}_{t}^{\top}\bm{J}_{t-1}\Delta \bm{x}_{t}}.
\end{equation}
In particular, the update $\bm{J}_{t}$ from \eqref{eq:BFGS} minimizes $\bm{J}_{t} \mapsto \| \bm{J}_{t-1} - \bm{J}_{t} \|_{F}$, subject to the secant condition \eqref{eq:secant-condition} and $\bm{J}_{t}$ being symmetric.  The second condition is because when $F=\bm{v} \mapsto \nabla f(\bm{v})$ for a function $f$,  then $\bm{J}_{\bm{x}} = \nabla^{2} f(\bm{x})$ is a symmetric matrix.  The inverse of $\bm{J}_{t}$ can be computed efficiently given the inverse of $\bm{J}_{t-1}$ following the Sherman-Morrison formula.

Under the same assumption as Theorem \ref{th:broyden-convergence} and the assumption that $\bm{J}_{\bm{x}^{*}}$ is symmetric definite positive,  the quasi-Newton method \eqref{alg:quasi-newton} with updates \eqref{eq:BFGS} converges locally and superlinearly to $\bm{x}^{*}$ (e.g., Theorem 5.5 and Corollary 5.6 in \cite{broyden1973local}). Note that \cite{rodomanov2021rates} provide explicit superlinear convergence rates (i.e. , bounds on $\| \bm{x}_{t} - \bm{x}^{*} \|_{2}$) for a large class of quasi-Newton methods, including BFGS.
\paragraph{Anderson mixing.} 
Let $m \in \N$.
In Anderson mixing \citep{anderson1965iterative},  the iterates $\left(\bm{x}_{t}\right)_{t \geq 0}$ are still $\bm{x}_{t+1} = \bm{x}_{t} - \bm{J}_{t}^{-1}F(\bm{x}_{t})$, but the information about the last $m+1$ iterates $\bm{x}_{t}, ..., \bm{x}_{t-m}$ and $F(\bm{x}_{t}), ...,F(\bm{x}_{t-m})$ is used to update the approximation $\bm{J}_{t}$. In particular, the computation of $\bm{J}_{t}$ relies on the last $m+1$ secant conditions:
\begin{equation}\label{eq:multi-secant-conditions}
\bm{J}_{t}\Delta \bm{x}_{i}= \Delta F_{i}, \forall \; i \in \{t, ...,t-m+1\}.
\end{equation}
With the notation $\Delta \bm{\XX}_{t} = \left( \Delta \bm{x}_{t-m+1}, ...,\Delta \bm{x}_{t} \right) \in \R^{n \times m,}  \Delta \bm{\FF}_{t} = \left( \Delta F_{t-m+1}, ...,\Delta F_{t} \right) \in \R^{n \times m},$ Equation \eqref{eq:multi-secant-conditions} can be concisely written as $\bm{J}_{t} \Delta \bm{\XX}_{t} = \Delta \bm{\FF}_{t}.$
Minimizing $\bm{J}_{t}  \mapsto \| \bm{J}_{t} - \bm{I} \|_{F}^{2}$ subject to \eqref{eq:multi-secant-conditions} yields a closed-form solution as
\begin{equation}\label{eq:good-anderson-update}\tag{And-I}
\bm{J}_{t} = \bm{I} + \left(\Delta \bm{\FF}_{t} - \Delta \bm{\XX}_{t} \right) \left(\Delta \bm{\XX}_{t}^{\top}\Delta \bm{\XX}_{t}\right)^{-1}\Delta \bm{\XX}_{t}^{\top}.
\end{equation}
The update \eqref{eq:good-anderson-update} relies on the \textit{type-I} Broyden update.
The \textit{type-II} Broyden update, which estimates $\bm{G}_{t}$ as an approximation of the inverse of the Jacobian,  satisfies the multi-secant condition in inverse form: $\bm{G}_{t}\Delta \bm{\FF}_{t} = \Delta \bm{\XX}_{t}.$ Subject to minimizing $\bm{G}_{t} \mapsto \| \bm{G}_{t} - \bm{I} \|_{F}^{2}$, the closed-form iteration for Anderson acceleration based on type-II Broyden update becomes
\begin{equation}\label{eq:bad-anderson-update}\tag{And-II}
\bm{G}_{t} = \bm{I} + \left(\Delta \bm{\XX}_{t} - \Delta \bm{\FF}_{t} \right) \left(\Delta \bm{\FF}_{t}^{\top}\Delta \bm{\FF}_{t}\right)^{-1}\Delta \bm{\FF}_{t}^{\top},
\end{equation}
and the update on $\bm{x}_{t}$ simply becomes $
\bm{x}_{t+1} = \bm{x}_{t} - \bm{G}_{t}F(\bm{x}_{t}).$
Note that it can be proved that $\bm{J}_{t}$ and $\bm{G}_{t}$ are still rank one updates from the previous approximations $\bm{J}_{t-1}$ and $\bm{G}_{t-1}$ (\cite{zhang2020globally},  Proposition 3.1, and  \cite{rohwedder2011analysis}, Theorem 3.2).
\paragraph{Superlinear Convergence of Anderson Acceleration.}
The relations between Anderson acceleration and Broyden's updates dates from \cite{eyert1996comparative}.
A proof of convergence of the full-memory version ($m=t$ at each iteration $t$) of \eqref{eq:good-anderson-update} can be found in \cite{gay1978solving}.  In particular, under the same assumption as Theorem \ref{th:broyden-convergence}, and with $m=t$ at each iteration $t$, if $\bm{x}_{0}$ is close enough to $\bm{x}^{*}$ then Anderson acceleration converges superlinearly to $\bm{x}^{*}$.
Local superlinear convergence, for a limited memory version ($m \in \N$ is fixed across iterations) of \eqref{eq:bad-anderson-update} is obtained in \cite{rohwedder2011analysis} (Theorem 5.2 and Theorem 5.4) under the same conditions.
Under a contraction assumption, \cite{toth2015convergence} obtain local linear convergence of \eqref{eq:bad-anderson-update} (Theorem 2.3 and Theorem 2.4). Crucially,  all these convergence results rely on the operator $F$ being differentiable, which is \textit{not} the case for MDPs (see the Bellman operator \eqref{eq:T_max-def}).  \cite{zhang2020globally} obtain convergence results for \eqref{eq:good-anderson-update} without the operator $F$ being differentiable, and this result is discussed more at length in Section \ref{sec:anderson-mdp}.
\end{document}